\documentclass[final,11pt]{elsarticle}

\usepackage[utf8]{inputenc} 

\usepackage{graphicx} 

\usepackage{CJKutf8}
\usepackage{tikz}
\usetikzlibrary{arrows,shapes,chains}
\usepackage{booktabs} 
\usepackage{array} 
\usepackage{paralist} 
\usepackage{verbatim} 

\usepackage{float}
\usepackage{epstopdf}
\usepackage{xfrac}
\usepackage{algorithm}
\usepackage{algorithmicx}
\usepackage{algpseudocode}
\usepackage{amsmath,dsfont}
\usepackage{harpoon}
\usepackage{amssymb}
\usepackage{mathrsfs}
\usepackage{float}
\usepackage{graphics}
\usepackage{color}
\usepackage{geometry}
\usepackage{epsfig}
\usepackage{diagbox}
\usepackage{mathtools}
\usepackage{comment}
\usepackage{bbm}
\usepackage{lineno}
\usepackage{fancyhdr}
\usepackage{graphicx}
\usepackage[draft]{changes}

\def\a{\alpha}
\def\b{\beta}

\def\dd{\mbox{d}}

\def\g{\gamma}

\def\m{\mathcal}

\def\x{{\bf x}}

\def\o{\omega}

\def\f{\frac}
\def\p{\partial}

\def\o{\omega}

\usepackage{wrapfig}
\usepackage{dsfont}
\usepackage{enumerate}
\usepackage{subfigure}
\usepackage{xcolor}
\usepackage{setspace}

\newtheorem{example}{\textup{\textbf{Example}}}
\begin{document}

\tikzstyle{format}=[rectangle,draw,thin,fill=white]
\tikzstyle{test}=[diamond,aspect=2,draw,thin]
\tikzstyle{point}=[coordinate,on grid,]

\allowdisplaybreaks
\DeclareGraphicsExtensions{.pdf,.jpeg,.png}
\numberwithin{equation}{section}
\begin{frontmatter}
\title{A frequency-dependent $p$-adaptive technique for spectral methods}

\author[UCLA]{Mingtao Xia}
\ead{xiamingtao97@g.ucla.edu}

\author[PKU]{Sihong Shao\corref{cor1}}
\ead{sihong@math.pku.edu.cn}

\author[UCLA]{Tom Chou\corref{cor1}}
\ead{tomchou@ucla.edu}

\cortext[cor1]{Corresponding author.}

\address[UCLA]{Department of Mathematics, UCLA, Los Angeles, CA
  90095-1555, USA}
\address[PKU]{LMAM and School of Mathematical
  Sciences, Peking University, Beijing 100871, CHINA}


\begin{abstract}
When using spectral methods, a question arises as how to determine the
expansion order, especially for time-dependent problems in which
emerging oscillations may require adjusting the expansion order.  In
this paper, we propose a frequency-dependent $p$-adaptive technique
that adaptively adjusts the expansion order based on a frequency
indicator.  Using this $p$-adaptive technique, combined with recently
proposed scaling and moving techniques, we are able to devise an
adaptive spectral method in unbounded domains that can capture and
handle diffusion, advection, and oscillations.  As an application, we
use this adaptive spectral method to numerically solve the
Schr\"{o}dinger equation in the whole domain and successfully capture
the solution's oscillatory behavior at infinity.
\end{abstract}
\end{frontmatter}

\section{Introduction}
\label{sec:intro}

Unbounded domain problems arise in many scientific applications
\cite{REVIEW,ADDER2020} and adaptive numerical methods are needed on
many occasions, for instance, in solving the Schr\"{o}dinger equation
in unbounded domains when the solution's behavior varies over time and
we wish to capture the solution's behavior in the whole domain.  As an
important class of numerical algorithms, adaptive methods have
witnessed numerous advances in their efficiency and accuracy
\cite{tang2003adaptive,babuska2012modeling,ren2000iterative,li2002adaptive}.
However, despite considerable progress that has been made for spectral
methods in unbounded domains \cite{shen2009some}, there are few
adaptive methods that apply in unbounded domains.

In \cite{Xia2020b}, adaptive scaling and moving techniques were
proposed for spectral methods in unbounded domains and it was noted
that adjusting the expansion order is necessary when the function
displays oscillatory behavior that varies over time.  In this paper,
we first develop a frequency-dependent technique for spectral methods
which adjusts the expansion order $N$ ($N+1$ basis functions are used
to approximate the solution). This technique takes advantage of the
frequency indicator defined in \cite{Xia2020b} and corresponds to
$p$-adaptivity \cite{Shao2011Adaptive,
  dumbser2007arbitrary,karihaloo2001accurate}. By adjusting the
expansion order efficiently, our $p$-adaptive technique can be used to
accurately solve problems with varying oscillatory.

By combining this $p$-adaptive technique with scaling and moving
methods, we develop an adaptive spectral method that can capture
diffusion, advection, and oscillations in unbounded domains.  Since
scaling and adjusting the expansion order both depend on the frequency
indicator, we also investigate the interdependence of these two
techniques. We demonstrate that appropriately adjusting the expansion
order can facilitate scaling to more efficiently distribute allocation points. In turn, proper
scaling can help avoid unnecessary increases in the expansion order when
it does not increase accuracy, thereby avoiding unnecessary
computational burden.

The significance of this adaptive spectral method is that it can
capture the solution's behavior in the whole domain. We demonstrate
the utility of our method by solving Schr\"{o}dinger's equation in
$\mathbb{R}$.  Here, the unboundedness and the oscillatory nature of
the solution pose two major numerical challenges
\cite{yang2014computation}.  Specifically, in the semiclassical
regime, when the wavelength of the solution is small, the function
becomes extremely oscillatory. Moreover, in certain situations, one
has to work with a very large computational domain that is difficult
to automatically determine.


Previous numerical methods which solve Schr\"{o}dinger's equation in
unbounded domains usually truncate the domain into a finite subdomain
and impose artificial boundary conditions, which may be nonlocal and
complicated \cite{han2005finite, yang2014computation, li2018stability,
  antoine2008review}.  Our adaptive spectral method tackles the
oscillatory problem directly in the original unbounded domain without
the need to truncate it or to devise an artificial boundary condition.


This paper is organized as follows. In the next section, we first
present a $p$-adaptive technique for spectral methods and use examples
to illustrate its efficiency. In Section \ref{unboundeddomain}, we
incorporate and study this technique within existing scaling and
moving techniques and devise an adaptive spectral method in unbounded
domains. Application of our adaptive spectral methods to numerically
solving Schr\"{o}dinger's equation is given in Section
\ref{schrodinger}. We summarize our results in Section
\ref{conclusion} and propose directions for future work.

\section{Frequency-dependent $p$-adaptivity}
\label{Boundeddomain}
We present a frequency-dependent $p$-adaptive spectral method based on
information extracted from only the numerical solution of
time-dependent problems. In \cite{Xia2020b}, we showed that a
frequency indicator defined for spectral methods is particularly
useful in measuring the contribution of high frequency modes in the
numerical solution. Because high frequency modes decay more slowly,
this indicator could be used to determine scaling in spectral methods
applied to unbounded domains. In this work, we will show that the frequency indicator can also be used
to determine whether more or fewer basis functions are needed to
refine or coarsen the numerical solution.

Given a set of orthogonal basis functions $\{B_{i}(x)\}_{i=0}^{\infty}$
under a specific weight function $\omega(x)>0$ in a domain $\Lambda$,
the frequency indicator associated with the interpolation of a
function 
\begin{equation}
\mathcal{I}_Nu(x)=U_N(x)=\sum_{i=0}^Nu_{i}B_{i}(x)
\end{equation}
is defined as in 
\cite{Xia2020b}
\begin{equation}
\m{F}(U_N) \coloneqq\left({\f{\sum\limits_{i=N-M+1}^{N}
\g_{i} u_{i}^2}{\sum\limits_{i=0}^{N}
\g_{i}u_{i}^2}}\right)^{\f{1}{2}},
\label{FREQSCAL}
\end{equation}
where $\gamma_{i} =\int_{\Lambda}B_{i}^2(x)\omega(x) \dd{x}$ is the
square of $L^2_{\omega}$-weighted norm of the basis function
$B_i(x)$. This frequency indicator measures the contribution of the $M$
highest-frequency components to the $L_{\o}^2$-weighted norm of
$U_N$. Here $M$ is often chosen to be $[\frac{N}{3}]$ following the
$\frac{2}{3}$-rule \cite{Hou2007Computing,Orszag1971On}.  This
indicator provides a lower bound for the error divided by the norm of
the numerical solution $\f{\|u -\m{I}_{N-M}u\|_{\o}}
{\|\m{I}_{N}u\|_{\o}}$ which is illustrated in \cite{Xia2020b}.  Thus,
the quality of the numerical interpolation $U_N$ can be measured by
$\m{F}(U_N)$.

For a time-dependent problem, the expansion order $N$ may need
adjusting dynamically, which can be reflected by the frequency
indicator.  If the frequency indicator increases, the lower bound for
$\f{\|u -\m{I}_{N-M}u\|_{\o}} {\|\m{I}_{N}u\|_{\o}}$ will also
increase.  On the other hand, as $N$ increases, the error $\|u
-\m{I}_{N}u\|_{\o}$ as well as $\m{F}(U_N)$ are expected to
decrease. By sufficiently increasing the expansion order $N$, the
frequency indicator as well as the error can be kept small.  If the
frequency indicator decreases, we can also consider decreasing $N$ to
relieve computational cost without compromising accuracy, as was done
in \cite{Shao2011Adaptive}.  The pseudo-code of the proposed
$p$-adaptive technique is given in Alg.~\ref{algrefining}.

\begin{algorithm}
\caption{\small Pseudo-code of the $p$-adaptive technique which may
  increase (refine) or decrease (coarsen) the expansion order $N$.}
\begin{algorithmic}[1]
\State Initialize $N, N_0$, $\gamma\geq1, \eta_0=\eta>1$, $\Delta t$, $T$, $\alpha$, $\beta$, $U_{N}(0)$, $N_{\max}$, $N_{\min}$
\State $t \gets 0$
\State $f_0 \gets \Call{frequency\_indicator}{U_{N}(t)}$ \label{alg1:i00}
\While{$t<T$}
\State $U_{N}(t+\Delta t)\gets \Call{evolve}{U_{N}(t),\Delta t}$\label{evolve}
\State $f\gets \Call{frequency\_indicator}{U_{N}(t+\Delta t)}$ \label{alg1:i0}
\State $l\gets0$
\If{$f > \eta f_0$} \quad\# refinement is needed \label{alg1:cond}
\While{$f > \eta f_0~ \textbf{and}  ~l\leq N_{\max} $} \label{alg1:whilecond}
\State $l\gets l+1$
\State $U_{N+1}  \gets  \Call{refine}{U_{N}(t+\Delta t)}$ \label{alg:refine0}
\State $N\gets N+1$
\State $f \gets \Call{frequency\_indicator}{U_N}$ \label{alg1:i1}
\EndWhile
\State $f_0\gets f$
\State $\eta\gets\gamma\eta$ \quad\# renew $\eta$
\ElsIf{$f<f_0/\eta_0$} \quad\# coarsening could be considered \label{whilereduce}
\State $r\gets\textbf{False}$
\While{$f<f_0/\eta_0~ \textbf{and}  ~N> N_{\min} ~\textbf{and not}~ r_1 $} 
\State $\tilde{U}_{N-1}(t+\Delta{t})\gets\Call{coarsen}{U_{N}(t+\Delta t)}$\label{alg:coarsen}
\State $f\gets\Call{frequency\_indicator}{\tilde{U}_{N-1}(t+\Delta t)}$
\If{$f<f_0$} \label{alg1:ifderefine}
\State $f_1\gets f$
\State $r\gets\textbf{True}$
\State $U_{N-1}(t+\Delta{t})\gets\tilde{U}_{N-1}(t+\Delta{t})$
\State $N\gets N-1$
\EndIf
\EndWhile
\If{$r$}
\State $f_0\gets f_1$
\EndIf
\EndIf
\State $t\gets t+\Delta t$
\EndWhile
\end{algorithmic}\label{algrefining}
\end{algorithm}

The $p$-adaptive spectral method in Alg.~\ref{algrefining} for
time-dependent problems consists of two ingredients: refinement
(increasing $N$) and coarsening (decreasing $N$).  It
maintains accuracy when there are emerging oscillations by increasing
the expansion order $N$. It also decreases $N$ when the expansion
order is larger than needed to avoid unnecessary computation.  In
Alg.~\ref{algrefining}, the $\Call{frequency\_indicator}{}$ subroutine
is to calculate the frequency indicator defined in
Eq.~\eqref{FREQSCAL} for the numerical solution $U_N$ while the
$\Call{evolve}{}$ subroutine is to obtain the numerical solution
$U_N(t+\Delta{t})$ at the next timestep from $U_N(t)$.

In Line \ref{alg:refine0} of Alg.~\ref{algrefining}, the
$\Call{refine}{}$ subroutine uses $U_N$ to generate a new numerical
solution with a larger expansion order $U_{N+1}$ (refine), and in Line
\ref{alg:coarsen} the $\Call{coarsen}{}$ subroutine uses $U_N$ to
generate a new numerical solution with a smaller expansion order
$U_{N-1}$ (coarsen). The refinement or coarsening is achieved by
reconstructing the function values of $U_{N+1}$ or $U_{N-1}$ at the
new set of collocation points $\{x_i\}$:

\begin{equation}
U_{N\pm1}(x_{i}, t)=U_N(x_{i}, t),\quad  i=0,...,N\pm 1, 
\end{equation}
where $U_{N+1}$ uses $N+2$ basis functions for refinement and
$U_{N-1}$ uses $N$ basis functions for coarsening.

In Alg.~\ref{algrefining}, $\eta f_0$ is the refinement threshold such
that if the current frequency indicator $f>\eta f_0$, we increase the
expansion order $N$. The \textbf{while} loop starting in Line
\ref{alg1:whilecond} ensures we either refine enough such that the
frequency indicator, after increasing $N$, is smaller than the
threshold $\eta f_0$, or the maximal allowable expansion order
increment within a single step $N_{\rm max}$ is reached.

After increasing $N$, $f_0$ is renewed to be the current frequency
indicator and $\eta$ is multiplied by a factor $\gamma\geq 1$,
enabling us to dynamically adjust the refinement threshold for the next
refinement in order to prevent increasing $N$ too fast without
substantially increasing accuracy.  On the other hand, when an
extremely large $N$ is needed to match the increasingly oscillatory
behavior of the numerical solution, we can set $\gamma\succeq 1$ or
even $\gamma=1$, as we will do in Examples~\ref{ex:simulatebounded} and
\ref{ex:schrodinger}.  We have observed numerically, as expected, that
the larger $\eta_0, \gamma$ are, the more difficult it is to increase
the expansion order.

We also consider reducing $N$ when a large expansion order is not
really needed and $f_0/\eta_0$ is the threshold for decreasing the
expansion order. If the condition in Line \ref{whilereduce} is
satisfied and $N>N_{\min}$, the minimal allowable expansion order, and
we have not increased $N$ in the current step, on the contrary we
consider decreasing the expansion order below Line \ref{whilereduce}.
As long as the frequency indicator of the new numerical solution with
the decreased expansion order $\m{F}(U_{N-1})$ is smaller than $f_0$,
the frequency indicator recorded after previously adjusting the
expansion order, reducing the expansion order is accepted; else
reducing the expansion order is declined. Therefore, $f_0$ after coarsening will not
surpass $f_0$ before coarsening. This procedure is described by the
\textbf{If} condition in Line \ref{alg1:ifderefine}.  If $N$ is
decreased, $f_0$ will also get renewed to be the latest frequency
indicator. 

In addition, if the current frequency indicator $f\in[\f{f_0}{\eta_0},
  \eta{f_0}]$, neither the refinement nor the coarsening subroutine
is activated.

Alg.~\ref{algrefining} can be generalized to higher dimensions in a
dimension-by-dimension manner.  The expansion order for each dimension
can change simultaneously within each timestep by using the tensor
product of one-dimensional basis functions, in much the same way
moving and scaling algorithms were generalized to higher dimensions
\cite{Xia2020b}.  For example, for a two-dimensional problem, given

\begin{equation}
U_{\vec{N}}(x, y):=\sum\limits_{i=0}^{N_x}\sum\limits_{j=0}^{N_y}u_{i, j}B_{i}(x)B_{j}(y)
\end{equation}
where $\vec{N}=(N_x, N_y)$, the frequency indicator in the
$x$-direction is defined as

\begin{equation}
\begin{aligned}
\mathcal{F}_x(U_{\vec{N}}) := \left(\f{\sum\limits_{i=N_x-M_x+1}^{N_x}
  \sum\limits_{j=0}^{N_y}\g_{i}
  \g_{j}u_{i, j}^2}{\sum\limits_{i=0}^{N_x}\sum\limits_{j=0}^{N_y}
\g_{i}\g_ju_{i, j}^2}\right)^{\f{1}{2}},
\end{aligned}
\label{freqdimension}
\end{equation}
while the frequency indicator in $y$-direction is similarly defined.
At each timestep, we keep $N_y$ fixed and use $\mathcal{F}_x$ to judge
whether or not to renew $N_x \to \tilde{N}_x$; simultaneously, we fix
$N_x$ and use $\mathcal{F}_y$ to renew $N_y \to \tilde{N}_y$ if
adjusting the expansion order in $y$ dimension is needed. Finally
$N_x, N_y$ are updated to $\tilde{N}_x, \tilde{N}_y$.

In this work, the relative $L^2_{\o}$-error
\begin{equation}
\text{Error} = \f{\|U_N-u\|_{\o}}{\|u\|_{\o}},
\label{RelaE}
\end{equation}
is used to measure the quality of the spectral approximation $U_N(x)$
compared to the reference solution $u(x)$.  Table \ref{tab:basis
  functions} lists some typical choices of orthogonal basis functions
for different domains $\Lambda$ that we use in this paper.

\begin{table}
  \centering
  \caption{\small Typical choices of basis functions
    $\{B_{i}\}_{i=0}^{\infty}$ and computational domain $\Lambda$.}
  \label{tab:basis functions}
   \newsavebox{\tablebox}
 \begin{lrbox}{\tablebox}
    \begin{tabular}{|c|c|c|c|c|}
\hline Computational domain & Bounded interval & $(0, \infty)$ &
  $(-\infty, \infty)$ \\ \hline Basis functions & Jacobi
  polynomials&Laguerre polynomials/functions& Hermite
  polynomials/functions\\ \hline
    \end{tabular}%
   \end{lrbox}
\scalebox{0.73}{\usebox{\tablebox}}
\end{table}

We provide two examples of using this $p$-adaptive technique in
Alg.~\ref{algrefining} below, where the generalized Jacobi polynomials
\cite{Spectral2011} are used.  Theorem 3.41 in \cite{Spectral2011}
gives an estimation for the interpolation error of a function $u$ in
the Jacobi-weighted Sobolev space for $\a, \b>-1$ as follows
\begin{equation}
\begin{aligned}
&\|\p_x^l(I_{N, \a, \b}u - u)\|_{\o_{\a+l,\b+l}}\leq \\
&\quad\:\:\quad\quad\quad\quad\quad\quad c\sqrt{\f{(N-m+1)!}{N!}}N^{l - (m+1)/2}\|\p_x^mu\|_{\o_{\a+m, \b+m}},\,\,
0\leq{l}\leq{m}\leq{N+1}, 
\end{aligned}
\end{equation}
where $c$ is a positive constant independent of $m, N$ and $u$.  When
$m>0$ and $l=0$, the left hand side becomes the interpolation error
$\|(I_{N, \a, \b}u - u)\|_{\o_{\a,\b}}$ which decreases with
$N$. Therefore, by increasing the expansion order for the Jacobi
polynomials it is generally true that the interpolation will be more
accurate. Theorem 7.16 and Theorem 7.17 in \cite{Spectral2011} give
similar error estimates for Laguerre and Hermite interpolations, which
reveals that under some smoothness assumptions, the interpolation
error decreases when the expansion order $N$ increases.


Since unbounded domain problems may involve diffusive and advective
behavior, we discuss and develop adaptive spectral methods in
unbounded domains in the next section.

 \begin{figure}
    \begin{center}
    \includegraphics[width=5.2in]{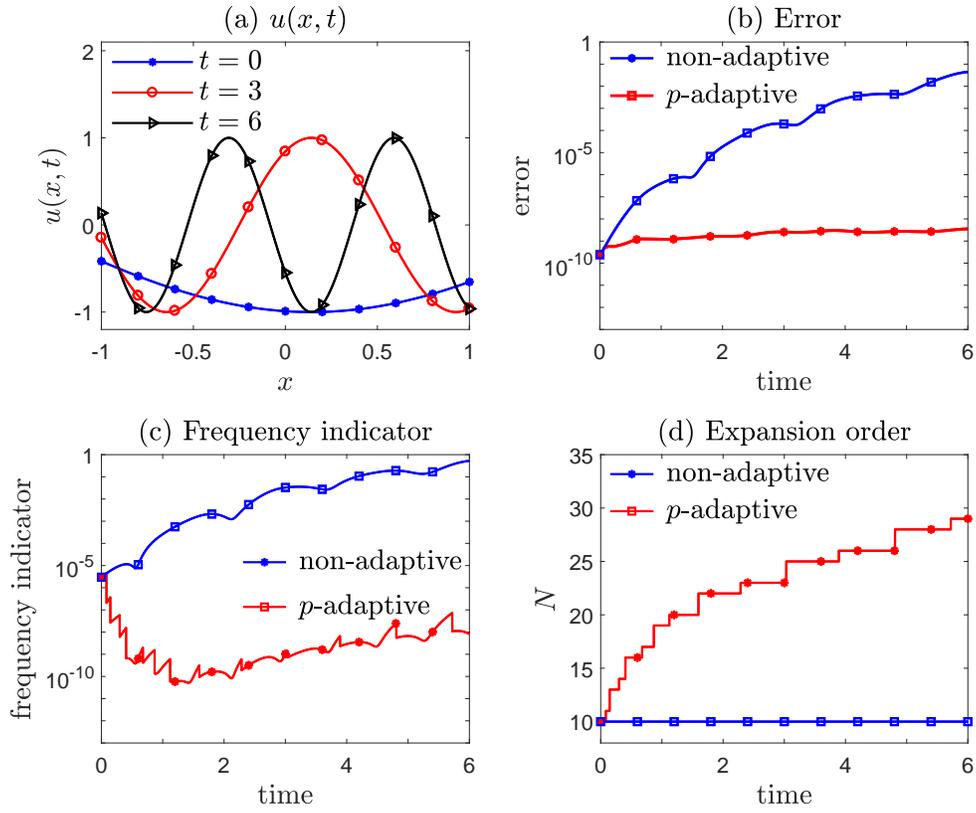}
\end{center}
\vspace{-3mm}
        \caption{\small Numerically solving Eq.~\eqref{chebyshev} with
          Chebyshev polynomials using Alg.~\ref{algrefining}. For
          solutions that become increasingly oscillatory, the
          $p$-adaptive technique can increase the expansion order
          effectively to capture the oscillations and maintain a small
          error by keeping the frequency indicator low. Using a fixed
          $N$ fails to maintain the frequency indicator and results in
          a large error.}
     \label{fig1}
\end{figure}

\begin{example}
\label{ex:Chebyshev}

\rm We numerically solve the PDE
\begin{equation}
\p_t u = \left(\f{x+2}{t+1}\right)\p_x u, \quad x\in[-1, 1],
\label{chebyshev}
\end{equation} 
with a Dirichlet boundary condition specified at $x=1$ given as $u(1,
t) = \cos 3(t+1)$. This PDE admits an analytical solution 
\begin{equation}
u(x, t) =\cos((t+1)(x+2)). 
\label{analytic}
\end{equation}
We solve it numerically by using Chebyshev polynomials with
Chebyshev-Gauss-Robatto quadrature nodes and weights.  The Chebyshev
polynomials are orthogonal under the weight function $\o(x) =
(1-x^2)^{-\f{1}{2}}$, \textit{i.e.}, they correspond to Jacobi
polynomials with $\a=\b=-\f{1}{2}$.  Since $u(x, t)$ becomes
increasingly oscillatory over time, an increasing expansion order is
required to capture these oscillations. We start with $N=10$ at $t=0$,
the parameters $\eta=1.5, \gamma=1.1, N_{\max}=3, N_{\min}=0$, and a
timestep $\Delta t=0.001$. We use a third order explicit Runge-Kutta
scheme to advance time.

The reference solution $u(x, t)$ is plotted in Fig.~\ref{fig1}(a).
The increasing oscillations lead to a fast rise in the frequency
indicator as the contribution from high frequency modes increases.
Keeping the same number of basis functions over time will fail as it
will be eventually incapable of capturing the shorter wavelength
oscillations.

However, a much more accurate approximation can be obtained (see
Fig.~\ref{fig1}(b)) with our $p$-adaptive method which maintains the
frequency indicator (see Fig.~\ref{fig1}(c)) by increasing the number
of basis functions (shown in Fig.~\ref{fig1}(d)). Furthermore, the
coarsening subroutine for decreasing the expansion order described in
the \textbf{while} loop in Line \ref{whilereduce} will not be
triggered (shown in Fig.~\ref{fig1}(d)).

When directly approximating the reference solution in
Eq.~\eqref{analytic}, we can achieve $10^{-8}$ accuracy with only 20
basis functions. However, when numerically solving
Eq.~\eqref{chebyshev}, the error will accumulate due to the increasing
oscillatory behavior which will require even more basis functions to
achieve the same accuracy as the direct approximation to
Eq.~\eqref{analytic}.  Thus, the oscillatory behavior of the solution
poses additional difficulties and requires even more refinement when
numerically solving a PDE.
\end{example}

\noindent Next, we present an example of a two-dimensional problem in $[-1,
  1]^2$.
 \begin{figure}
\begin{center}
      \includegraphics[width=5.2in]{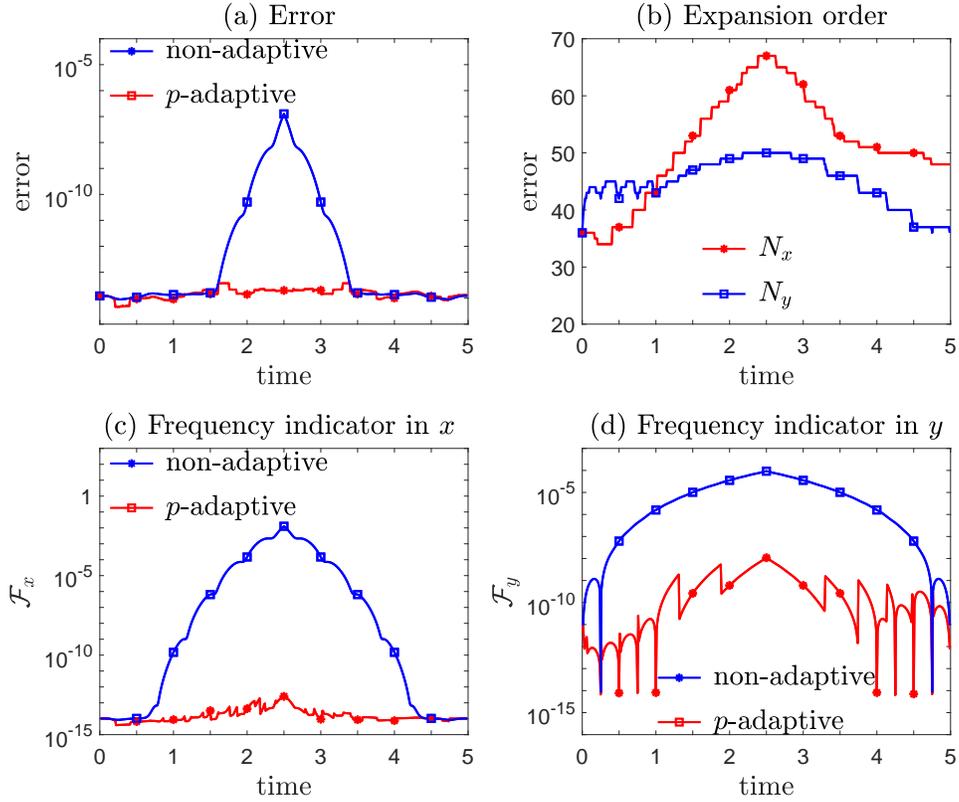}
\end{center}
\vspace{-3mm}
        \caption{\small Using the $p$-adaptive technique to
          approximate the two-dimensional function in
          Eq.~\eqref{2dapprox} with Legendre polynomials.  Refinement
          is applied in each direction simultaneously to capture
          increasing oscillations in both directions. Coarsening is
          applied when large expansion orders are not
          needed. Anisotropic oscillatory behavior requires adjusting
          the expansion order in each direction differently. The
          frequency indicators in both dimensions are kept low,
          leading to a small error.}
     \label{fig2}
\end{figure}

\begin{example}
\label{ex:Lagendre2D}
\rm
We approximate the function
\begin{equation}
u(x, y, t) =
\cos\left(xy(5-2\vert t-\tfrac{5}{2}\vert)\right) + y^{10 - 4|t-5/2|}
\sin \left(4x(5-2\vert t-\tfrac{5}{2}\vert)\right), \,\,  (x, y)\in [-1, 1]^2
\label{2dapprox}
\end{equation}
by Legendre polynomials (corresponding to Jacobi polynomials with
$\a=\b=0$) with Legendre-Gauss-Robatto quadrature nodes and weights in
both dimensions.  Within $t\in[0, \f{5}{2}]$, the function becomes
more oscillatory over time in both dimensions, requiring increasing
expansion orders.  For $t\in[\f{5}{2}, \f{7}{2}]$, the error for approximation with fixed expansion orders in both dimensions
decreases because the function becomes less
oscillatory, and therefore a reduction in expansion orders in both directions can
be used to reduce computational effort without compromising
accuracy. Since the function is not symmetric in $x$ and $y$ and the
adjustment of expansion is anisotropic. We show that
Alg.~\ref{algrefining} can appropriately increase $N_x, N_y$ when
$t<\f{5}{2}$ and reduce $N_x, N_y$ when $t\geq\f{5}{2}$.  We take
$N_x=N_y=36$ at $t=0$ with a timestep $\Delta{t}=0.01$, and
$\gamma_x=\gamma_y=1.1, \eta_x=\eta_y=1.1, N_{\max, x}=N_{\max, y}=3,
N_{x,\min}=N_{y,\min}=0$.

It is clear from Fig.~\ref{fig2}(a) that fixing the number of basis
functions in each dimension leads to an approximation which
deteriorates while the proposed $p$-adaptive spectral method is able to keep the
error small. Furthermore, when $t\in[\f{5}{2}, 5]$ we see that with
fixed expansion order the approximation error decreases, indicating
that coarsening may be performed to relieve computational burden while
maintaining accuracy.  Alg.~\ref{algrefining} first tracks the
increasing oscillation by increasing expansion orders in both $x$ and
$y$ dimensions. When $t\geq\f{5}{2}$, Alg.~\ref{algrefining} senses a
decrease in frequency indicator and decreases both $N_x$ and $N_y$
adaptively (shown in Figs.~\ref{fig2}(b)) without compromising
accuracy, as shown in Fig.~\ref{fig2}(a).

Since $\sin(4x(5-2|t-\f{5}{2}|))$ is the most oscillatory term in
$u(x, y, t)$, the function becomes more oscillatory in $x$ than in $y$
when $t\in[0, \f{5}{2}]$.  Because the function displays different
oscillatory behavior in $x$- and $y$-directions, the expansion order
should be adjusted anisotropically and $N_x$ needs increasing more
than $N_y$ in order to maintain $\mathcal{F}_x$ small.  Both
$\mathcal{F}_x$ and $\mathcal{F}_y$ are maintained well for the
$p$-adaptive approximation over time (shown in Fig.~\ref{fig2}(c) and (d)),
leading to satisfactory error control. Overall, Alg.~\ref{algrefining}
preserves accuracy for all times while still avoids using excessive
values of $N_x$ and $N_y$ when they are not needed.

\end{example}

\section{Adaptive spectral methods in unbounded domains}
\label{unboundeddomain}
\begin{figure}  
\centering
		\scriptsize  
		\tikzstyle{format}=[rectangle,draw,thin,fill=white]  
		\tikzstyle{test}=[diamond,aspect=2,draw,thin]  
		\tikzstyle{point}=[coordinate,on grid,]  
		\begin{tikzpicture}
		\node[format] (start){Initialize $N$, $\Delta t$, $T$, $\beta$, $U_{N}^{(\beta)}(0)$, $x_L, x_R = x^{(\beta)}_{[\f{N+2}{3}]}$};
		\node[format,below of=start,node distance=10mm] (define0){$f_0, f_1 \gets \Call{frequency\_indicator}{U_{N, x_L}^{(\beta)}(x, t)}$};
		\node[format,below of=define0,node distance=10mm] (define){$e_0 \gets \Call{exterior\_error\_indicator}{U_{N, x_L}^{(\b)}(x, 0), x_R}$};
		\node[test,below of=define,node distance=10mm](time){$t<T$};
		\node[format,left of=time,node distance=25mm](over){End};
		\node[test,below of=time,node distance=12mm] (PCFinit){MOVE?};
		\node[point, left of=PCFinit, node distance = 25mm](moveinter1){};
		\node[format,below of=moveinter1,node distance=5mm](movee0){Renew $e_0, x_L, x_R$};
		\node[point, below of=movee0, node distance=5mm](pointinter1){};
		\node[point, right of=pointinter1, node distance=25mm](pointinter2){};
		\node[test,below of=pointinter2,node distance=7mm] (DS18init){SCALE?};
		\node[point, left of=DS18init, node distance=25mm](scaleinter){};
		\node[format,below of=scaleinter,node distance=5mm] (LCDinit){Renew $\beta, f_1, x_R$};
		\node[point, below of=LCDinit, node distance=5mm](pointscale){};
		\node[point, right of=pointscale, node distance = 25mm](pointscale2){};
		\node[test,below of=pointscale2,node distance=13mm](setkeycheck){REFINE or COARSEN?};
		\node[point,right of=PCFinit,node distance=18mm](movepoint){};
		\node[point,below of=movepoint,node distance=10mm](movepoint2){};
		\node[point,left of=movepoint2,node distance=18mm](movepoint3){};
		\node[point,left of=setkeycheck,node distance=15mm](point3){};
		\node[point,right of=point3,node distance=15mm](point4){};
		\node[point, right of=setkeycheck, node distance=40mm](pointconnect1){};
		\node[format, above of=pointconnect1, node distance=20mm](increasetime){$t = t + \Delta{t}$};
	     \node[point, above of=increasetime, node distance=32mm](pointconnect2){};
		\node[test,below of=setkeycheck, node distance=18mm](processtime1){REFINE?};
		\node[point, left of=processtime1, node distance=25mm](refineinter){};
		\node[point, right of=DS18init, node distance=18mm](scalepoint1){};

		\node[point, below of=scalepoint1, node distance=10mm](scalepoint2){};
		\node[point, left of=scalepoint2, node distance=18mm](scalepoint3){};
		\node[point,right of=processtime1,node distance=20mm](refinepoint){};
		\node[point,below of=refinepoint,node distance=10mm](refinepoint2){};
		\node[point,left of=refinepoint2,node distance=15mm](refinepoint3){};
		\node[format,below of=refineinter, node distance=5mm](gettemp){Renew $\eta$ };
		\node[point, below of=gettemp, node distance=5mm](refinepoint4){};
		\node[point, right of=refinepoint4, node distance=20mm](refinepoint5){};
		\node[format,below of=processtime1](display){Renew $e_0, f_0, f_1, x_R, N$};
		\node[point,below of=display,node distance=5mm](point1){};
		\node[point, right of=point1, node distance=40mm](refineover){};
		\draw[->](start)--(define0);
		\draw[->](define0)--(define);
		\draw[->](define)--(time);
		\draw[-](movee0)--(pointinter1);
		\draw[-](LCDinit)--(pointscale);
		\draw[-](pointinter1)--(pointinter2);
		\draw[->](time)--node[left]{Yes}(PCFinit);
		\draw[->](time)--node[above]{No}(over);
		\draw[->](pointconnect1)--(increasetime);
		\draw[-](increasetime)--(pointconnect2);
		\draw[-](PCFinit)--node[above]{No}(movepoint);
		\draw[-](movepoint)--(movepoint2);
		\draw[->](moveinter1)->(movee0);
		\draw[-](DS18init)--node[above]{Yes}(scaleinter);
		\draw[->](scaleinter)->(LCDinit);
		\draw[-](PCFinit)--node[above]{Yes}(moveinter1);
		\draw[->](pointinter2)--(DS18init);
		\draw[-](movepoint2)--(pointinter2);
		\draw[-](DS18init)--node[above]{No}(scalepoint1);
		\draw[-](scalepoint1)--(scalepoint2);
		\draw[-](processtime1)--node[above]{No}(refinepoint);
		\draw[-](refinepoint)--(refinepoint2);
		\draw[-](processtime1)--node[above]{Yes}(refineinter);
		\draw[->](refinepoint2)--(display);
		\draw[-](scalepoint2)--(pointscale2);
		\draw[->](pointscale2)->(setkeycheck);
		\draw[-](pointscale)--(pointscale2);
		\draw[-](setkeycheck)--node[above]{No}(pointconnect1);
		\draw[->](pointconnect2)--(time);
		\draw[->](refineinter)--(gettemp);
		\draw[->](setkeycheck) --node[left]{Yes}(processtime1);
		\draw[-](display)--(point1);
		\draw[-](gettemp)--(refinepoint4);
		\draw[->](refinepoint4)--(display);
		\draw[-](point1)--(refineover);
		\draw[-](refineover)--(pointconnect1);
		\end{tikzpicture}  
\caption{\small Flow chart of an adaptive spectral method in unbounded
  domains which consists of moving, scaling, refinement and
  coarsening techniques.}
\label{algmovingscaling}
	\end{figure}
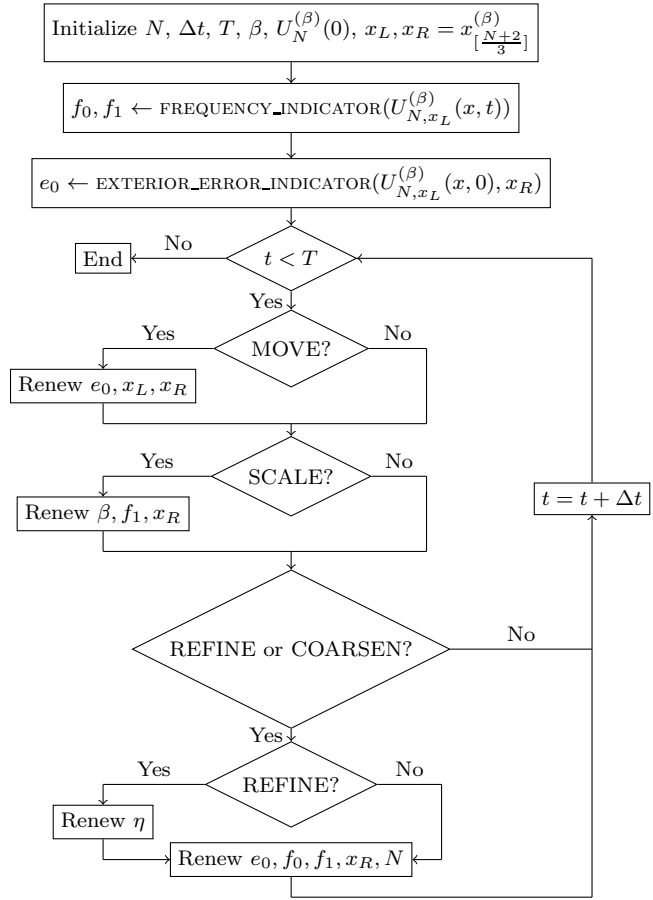

Unbounded domain problems are often more difficult to numerically
solve than bounded domain problems. Diffusion and advection in
unbounded domains necessitates knowledge of the solution's behavior at
infinity.  To distinguish and handle diffusive and advective behavior
in unbounded domains, techniques for scaling and moving basis
functions are proposed in \cite{Xia2020b}. When combining scaling,
moving, refinement and coarsening, we can devise a comprehensive
adaptive spectral approach for unbounded domains. A flow chart of our
overall approach is given in Fig.~\ref{algmovingscaling}.  The
scaling, refinement and coarsening techniques all rely on a common
frequency indicator.

As is stated in \cite{Xia2020b}, advection may cause a false increase
in the frequency indicator. Thus, we must first compensate for
advection by the moving technique before we consider either scaling or
adjusting the expansion order. Next, as the cost of changing the
scaling factor is lower than increasing the expansion order, we
implement scaling before adjusting the expansion order.  Only if
scaling cannot maintain the frequency indicator below the refinement
threshold do we consider increasing the expansion order. Coarsening is
also considered after scaling if the frequency indicator decreases
below the threshold for coarsening while no refinement is performed in
the current timestep.

As we have done in \cite{Xia2020b}, we also decrease the scaling
factor $\beta$ by multiplying it by a common ratio $q<1$ if the
current frequency indicator is larger than the scaling threshold
$f>\nu f_1$. The scaling we perform here contains an additional step:
when the current frequency indicator decreases and is below $f_1$, we
consider increasing the scaling factor $\beta$ by dividing it by the
common ratio $q$ as long as the frequency indicator decreases after
increasing $\beta$. When $f\in[f_1, \nu f_1]$, $\beta$ is neither
increased nor decreased. Thus, at each step, the scaling factor
$\beta$ may be either increased or decreased as long as the frequency
indicator decreases after adjusting the scaling factor.  A decrease in the
scaling factor indicates that the allocation points are more
efficiently distributed. These changes avoid unnecessary computational
burden that may arise if $N$ is excessively increased.  We briefly
describe our modified scaling subroutine for one timestep in
Alg.~\ref{algscaling}.
  
\begin{algorithm}
\caption{\small Pseudo-code of the frequency-dependent scaling
  technique which may increase or decrease the scaling factor
  $\beta$.}
\begin{algorithmic}[1]
\State $f\gets \Call{frequency\_indicator}{U_{N}^{(\alpha,\beta)}(t+\Delta t)}$ \label{alg:i0}
\If{$f > \nu f_1$} \label{alg:cond} \quad \# try decreasing $\beta$
\State $\tilde{\b}  \gets q \beta$ 
\State $U_N^{(\a,\tilde{\b})}  \gets  \Call{scale}{U_{N}^{(\alpha,\beta)}(t+\Delta t),\tilde{\b}}$ \label{alg:s0}
\State $\tilde{f} \gets \Call{frequency\_indicator}{U_N^{(\a,\tilde{\b})}}$ \label{alg:i1}
\While{$\tilde{f}\leq f$ \textbf{and} $\tilde{\b}\geq \underline{\b}$}\label{alg:scalewhile}
\State $\b \gets \tilde{\b}$
\State $U_{N}^{(\alpha,\beta)}(t+\Delta t)\gets U_N^{(\a,\tilde{\b})}$
\State $f_1 \gets \tilde{f}$
\State $f \gets \tilde{f}$
\State $\tilde{\b}  \gets q \beta$ \label{q}
\State $U_N^{(\a,\tilde{\b})}  \gets  \Call{scale}{U_{N}^{(\alpha,\beta)}(t+\Delta t),\tilde{\b}}$  \label{alg:s1}
\State $\tilde{f} \gets \Call{frequency\_indicator}{U_N^{(\a,\tilde{\b})}}$ \label{alg:i22}
\EndWhile
\ElsIf{$f<f_1$}\label{increasebeta} \quad \# try increasing $\beta$
\State $\tilde{\b}  \gets \beta/q$  
\State $U_N^{(\a,\tilde{\b})}  \gets  \Call{scale}{U_{N}^{(\alpha,\beta)}(t+\Delta t),\tilde{\b}}$ \label{alg:s02}
\State $\tilde{f} \gets \Call{frequency\_indicator}{U_N^{(\a,\tilde{\b})}}$ \label{alg:i12}
\While{$\tilde{f}\leq f$ \textbf{and} $\tilde{\b}\leq \overline{\b}$}\label{alg:scalewhile2}
\State $\b \gets \tilde{\b}$
\State $U_{N}^{(\alpha,\beta)}(t+\Delta t)\gets U_N^{(\a,\tilde{\b})}$
\State $f_1 \gets \tilde{f}$
\State $f \gets \tilde{f}$
\State $\tilde{\b}  \gets \beta/q$
\State $U_N^{(\a,\tilde{\b})}  \gets  \Call{scale}{U_{N}^{(\alpha,\beta)}(t+\Delta t),\tilde{\b}}$ 
\State $\tilde{f} \gets \Call{frequency\_indicator}{U_N^{(\a,\tilde{\b})}}$
\EndWhile
\EndIf
\end{algorithmic}\label{algscaling}
\end{algorithm}

For simplicity, we assume that the function is moving rightward so we
need to move the basis functions rightward. Therefore, $(x_R, \infty)$
is the ``exterior domain" of the spectral approximation on which we
wish to control the error as illustrated in \cite{Xia2020b}.  For
Laguerre polynomials/functions the parameter $x_L$ in the algorithm in
Fig.~\ref{algmovingscaling} denotes the starting point for the
approximation, while for Hermite polynomials/functions $x_L$
represents the translation of Hermite polynomials/functions,
\textit{i.e.}, we use $\{\mathcal{H}_i(x-x_L)\}$ or
$\{\hat{\mathcal{H}}_i(x-x_L)\}$.  After the expansion order $N$ has
changed, we need to renew both the threshold for scaling and the threshold for adjusting
the expansion order. We also need to renew the threshold for moving,
denoted by $e_0$, after $N$ has changed because different $N$ leads to
different $x_R$.  $U_{N, x_L}^{(\b)}$ is the spectral approximation
with the scaling factor $\beta$. Here the exterior-error indicator for
the semi-unbounded domain is defined in \cite{Xia2020b} and we can
generalize it to $\mathbb{R}$ when using Hermite
polynomials/functions

\begin{equation}
   \mathcal{E}(U_{N, x_L}^{(\b)},x_R) = \f{\|\p_x U_{N, x_L}^{(\b)}
\cdot\mathbb{I}_{(x_R,+\infty)}\|_{{\o}_{\beta}}}{\|\p_x U_{N, x_L}^{(\b)}\cdot
     \mathbb{I}_{(-\infty,+\infty)}\|_{{\o}_{\beta}}},
\label{errorindicator}
 \end{equation}
where $\o_{\beta}$ is the weight function and $x_R$ is taken to be
$x_{[\f{2N+2}{3}]}^{(\b)}$ for Hermite functions/polynomials and
$x_{[\f{N+2}{3}]}^{(\a, \b)}$ for Laguerre functions/polynomials
\cite{Xia2020b} in view of the often-used $\f{2}{3}$-rule. The
difference between the choices of $x_R$ for Hermite and Laguerre basis
functions arises because the allocation points for Hermite functions
are symmetrically distributed around their center while those for
Laguerre functions are one-sided, to the right of the starting point
$x_L$ in the axis.

For the scaling subroutine we need the following parameters: the
common ratio $q<1$ that we use to geometrically shrink/increase the
scaling factor, the parameter describing the threshold for considering
shrinking the scaling factor $\nu$; a predetermined lower bound for
the scaling factor $\underline{\beta}$ and an upper bound
$\overline{\beta}$. For the moving subroutine, required parameters
include the minimal displacement for the moving technique $\delta$,
the maximal displacement within a single timestep $d_{\rm max}$, and
the parameter of the threshold for activating the moving technique
$\mu$. The $\Call{scale}{}$ and $\Call{move}{}$ in
Fig.~\ref{algmovingscaling} are the scaling and moving subroutines and
$\Call{exterior\_error\_indicator}{}$ calculates the exterior-error
indicator for the moving subroutine. Detailed discussions about
scaling and moving techniques are given in \cite{Xia2020b}.

After first applying the moving technique, adjusting the expansion
order and scaling both depend on the frequency indicator and aim to
keep the frequency indicator low to control the error. The
relationship and interdependence between them is key to understanding
and justifying the first-scaling-then-adjusting expansion-order
procedure in Fig.~\ref{algmovingscaling}.  Thus, we need to
investigate how the proposed scaling technique will affect our
$p$-adaptive technique and how these two techniques interact with each
other. We use two examples containing both diffusive and oscillatory
behavior to investigate how the two techniques will be activated and
influence each other.  In Example \ref{ex:sensitivityscaling}, both
refinement and reducing $\beta$ are needed for matching increasing
oscillatory and diffusive behavior of the solution; in Example
\ref{ex:coarseningscaling}, a less oscillatory and diffusive solution
over time implies that coarsening and increasing $\beta$ may be
considered.

\begin{table}
  \centering
  \caption{\small Error, $\beta$, and $N$ at $t=5$ for different
    $\eta$ and $\gamma$ with both $p$-adaptive and scaling
    techniques.}
  \label{tab:refinescale}
   \begin{lrbox}{\tablebox}
    \begin{tabular}{|l|r|r|r|r|}
\hline
\diagbox{$\gamma$}{$\eta$} & 1.2 & 1.5 & 2& 4 \\
\hline
1.05 & \diagbox{1.305e-05}{1.434, 67}&\diagbox{2.346-05}{1.434, 64}&\diagbox{7.687e-05}{1.362, 58}&\diagbox{6.030e-05}{2.053, 55}\\
\hline
1.1&\diagbox{2.500e-05}{1.510, 62}&\diagbox{5.396e-05}{1.434, 69}  & \diagbox{5.513e-05}{1.673, 66} &\diagbox{6.030e-05}{2.053, 55} \\ 
\hline
1.2 & \diagbox{4.451e-05}{1.853, 56}&\diagbox{5.512e-05}{1.673, 56}   & \diagbox{8.927e-05}{1.673, 53}  &\diagbox{8.706e-05}{1.853, 53} \\
\hline
1.5 &\diagbox{1.369e-04}{1.589, 52}&\diagbox{8.927e-05}{1.673, 53}  & \diagbox{1.099e-04}{1.761, 52}  &\diagbox{8.702e-05}{1.951, 53} \\
\hline
    \end{tabular}
       \end{lrbox}
\scalebox{0.8}{\usebox{\tablebox}}
\end{table}

\begin{example}
\label{ex:sensitivityscaling}
\rm

We approximate the function
\begin{equation}
u(x, t) = \exp\left[-\f{x}{(bt+a)}\right]\cos x, \,\,\, t\in\mathbb{R}^+
\end{equation}
with the generalized Laguerre function basis
$\{\hat{\mathcal{L}}^{(\alpha, \beta)}_i(x)\}_{i=0}^{\infty}$ discussed
in \cite{Xia2020b} with the parameter $\alpha=0$.  The magnitude of
oscillations for this function, $\exp(-\f{x}{(bt+a)})$, increases over
time, requiring proper scaling. Under a variable transformation $y =
\frac{x}{bt+a}$, $u(x,t)$ can be rewritten as $u(y, t) =
\cos\left((bt+a)y\right)\exp(-y)$, indicating that the solution is
increasingly oscillatory in $y$ as time increases.  Thus, if we reduce
the scaling factor $\beta$ to match the diffusive behavior of the
solution, proper refinement is also required. In other words,
diffusive and oscillatory behavior is coupled in this example. We
carry out numerical experiments using the algorithm described in
Fig.~\ref{algmovingscaling} with different $(\eta, \gamma)$ to
investigate how scaling and refinement influence each other.  We
deactivate the moving technique by setting $d_{\rm max}=0$ since the
solution exhibits no intrinsic advection. Even if we had allowed
moving, it was hardly activated. We set $\Delta{t}=10^{-3}, N=50$ at
$t=0$ and $a=2, b=0.7$. $q=v^{-1}=0.95, \underline{\beta}=0.3,
\overline{\beta}=10, N_{\min}=0, N_{\max}=3$ and choose the initial
scaling factor $\beta=4$.

\begin{table}
  \centering
  \caption{\small Error and $N$ at $t=5$ for different $\eta$ and
    $\gamma$ with the $p$-adaptive technique but without the scaling
    technique, $\beta=4$.}
  \label{tab:refinenoscale}
   \begin{lrbox}{\tablebox}
    \begin{tabular}{|l|r|r|r|r|}
\hline
\diagbox{$\gamma$}{$\eta$} & 1.2 & 1.5 & 2& 4 \\
\hline
1.05 & \diagbox{1.316e-05}{79}&\diagbox{3.720e-05}{71}&\diagbox{9.724e-05}{65}&\diagbox{2.364e-04}{59}\\
\hline
1.1& \diagbox{4.372e-05}{70}&\diagbox{6.544e-05}{67}  & \diagbox{9.823e-05}{64} &\diagbox{2.607e-04}{58} \\ 
\hline
1.2 & \diagbox{9.724e-05}{65}&\diagbox{1.534e-04}{62}   & \diagbox{1.597e-04}{61}  &\diagbox{2.607e-04}{58} \\
\hline
1.5 &\diagbox{2.364e-03}{59} &\diagbox{2.364e-04}{59}  & \diagbox{2.607e-04}{58}  &\diagbox{3.508e-04}{56} \\
\hline
    \end{tabular}
           \end{lrbox}
\scalebox{0.8}{\usebox{\tablebox}}
\end{table}

In Table \ref{tab:refinescale} and Table \ref{tab:refinenoscale} the
error in $\mathbb{R}^+$ is recorded in the lower-left part of each
entry while the scaling factor $\beta$ and expansion order $N$ at
$t=5$ is recorded in the upper-right. By comparing entries in each
column/row for smaller $\eta,\gamma$, both tables show the expansion
order $N$ is likely to be increased more when the threshold for
refinement $\eta f_0$ is lower.

It can be observed from Table \ref{tab:refinescale} that with more
refinement $\beta$ tends to be smaller. This interaction between
$p$-adaptivity and scaling arises because more refinement leads to a
larger expansion order $N$ and a smaller scaling threshold $\nu f_1$.
Since scaling will only be performed if the frequency indicator after
scaling decreases, proper refinement is not likely to lead to
over-scaling.

Moreover, by comparing $N$ at $t=5$ between Tables
\ref{tab:refinescale} and \ref{tab:refinenoscale}, we see that $N$
tends to be smaller with the scaling technique for the same $\gamma,
\eta$. This implies that without scaling, the refinement procedure is
more often activated, leading to a larger $N$ to compensate for the
incapability of scaling alone to maintain a low frequency
indicator. This results in a larger computational burden without an
improvement in accuracy.  This behavior has been expected from the
design of Alg.~\ref{algmovingscaling} since we put scaling before
refinement so that redistribution of collocation points is tried first
to avoid unnecessary refinement when the increase in frequency
indicator results from diffusion instead of oscillation.

\end{example}

\begin{example}
\label{ex:coarseningscaling}
\rm

We approximate the function
\begin{equation}
u(x, t) = \exp\left[-(bt+a)x\right]\cos x,\,\,  x, t\in\mathbb{R}^+
\label{coarsescale}
\end{equation}
with the generalized Laguerre function basis with the parameter
$\alpha=0$.  The magnitude of oscillations for this function,
$\exp(-(bt+a)x)$, decreases over time and increasing the scaling
factor $\beta$ to more densely redistribute the allocation points is
needed. Furthermore, under the variable transformation $y = (bt+a)x$,
$u(x,t)$ can be rewritten as $u(y, t) = \cos(\f{y}{bt+a})\exp(-y)$.
Since the oscillations decrease with $y$, one can reduce the expansion
order. We consider coarsening with or without scaling to investigate
whether increasing $\beta$ can facilitate coarsening (and save
computational effort) or result in higher accuracy. We carry out
numerical experiments using the algorithm described in
Fig.~\ref{algmovingscaling} and different $(\eta, \gamma)$ and also
deactivate the moving technique by setting $d_{\rm max}=0$ since the
solution exhibits no intrinsic advection. We set $\Delta{t}=10^{-3},
N=50$ at the beginning and set the parameters $a=\f{1}{2}, b=0.5$,
$q=v^{-1}=0.95, \underline{\beta}=0.3, \overline{\beta}=10,
N_{\min}=0, N_{\max}=3$ and initial scaling factor $\beta=4$. We use a
different threshold $\eta_0$ for coarsening and we have checked
numerically that the parameter $\gamma$ in the refinement subroutine
will not affect the coarsening subroutine in this example.
\begin{table}
  \centering
  \caption{\small Error, $\beta$ and $N$ at $t=5$ for different
    $\eta_0$ and $\gamma$ with/without scaling for the $p$-adaptive
    technique.}
  \label{tab:coarsenscale}
     \begin{lrbox}{\tablebox}
    \begin{tabular}{|l|r|r|r|r|}
\hline
$\eta$& 1.2 & 1.5 & 2& 4 \\
\hline
Scaled & \diagbox{6.514e-10}{5.728, 11}&\diagbox{8.885e-12}{7.032, 13}&\diagbox{1.260e-11}{6.347, 13}&\diagbox{1.255e-14}{6.681, 17}\\
\hline
Unscaled &\diagbox{2.707e-12}{4, 20}&\diagbox{8.127e-11}{4, 28}  & \diagbox{9.417e-15}{4, 31} &\diagbox{9.672e-15}{4, 30} \\ 
\hline
    \end{tabular}
               \end{lrbox}
\scalebox{0.8}{\usebox{\tablebox}}
\end{table}

In Table \ref{tab:coarsenscale} the error in $\mathbb{R}^+$ is
recorded in the lower-left part of each entry while the scaling factor
$\beta$ and expansion order $N$ at $t=5$ is recorded in the
upper-right. By comparing entries in each row we see that a smaller
$\eta_0$ will lead to easier coarsening and a smaller $N$ at
$t=5$. Since the approximation with larger $N$ is always better,
whether we can achieve the same level of accuracy with a smaller
expansion order $N$ if proper scaling is implemented is of
interest. The initial approximation error is $1.960\times 10^{-9}$ and
the approximation will not worsen after coarsening regardless of
$\eta_0$ because in the $p$-adaptive subroutine coarsening is allowed
only when the post-coarsening frequency indicator remains below the
previous threshold $f_0$. Moreover, by comparing the two rows in Table
\ref{tab:coarsenscale} we see that if the solution concentrates and
becomes less diffusive, increasing $\beta$ and more efficiently
redistributing the allocation points allows the scaling technique to achieve high
accuracy with fewer expansion orders than without the scaling technique.
  \begin{figure}
\begin{center}
      \includegraphics[width=5in]{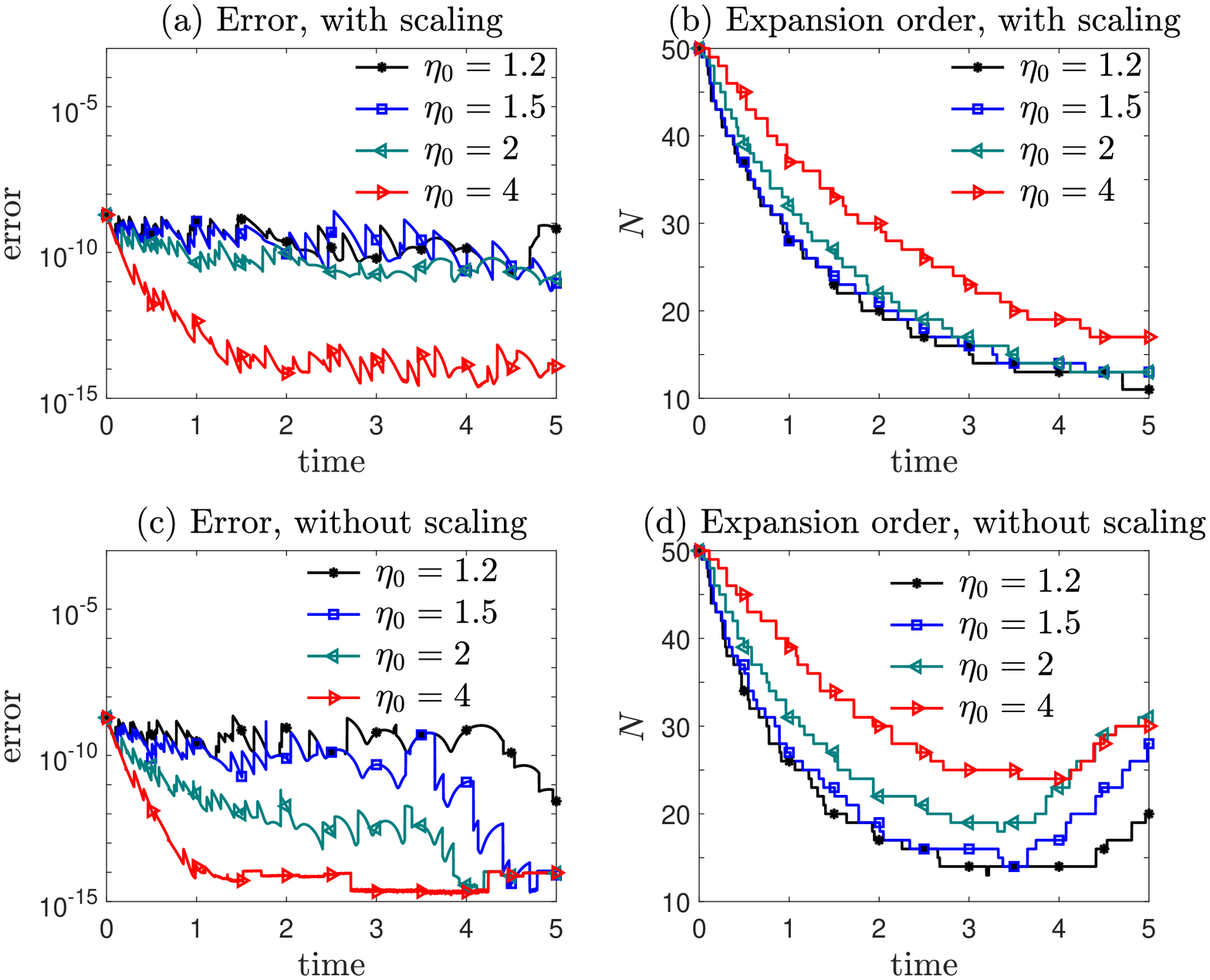}
\includegraphics[width=5in]{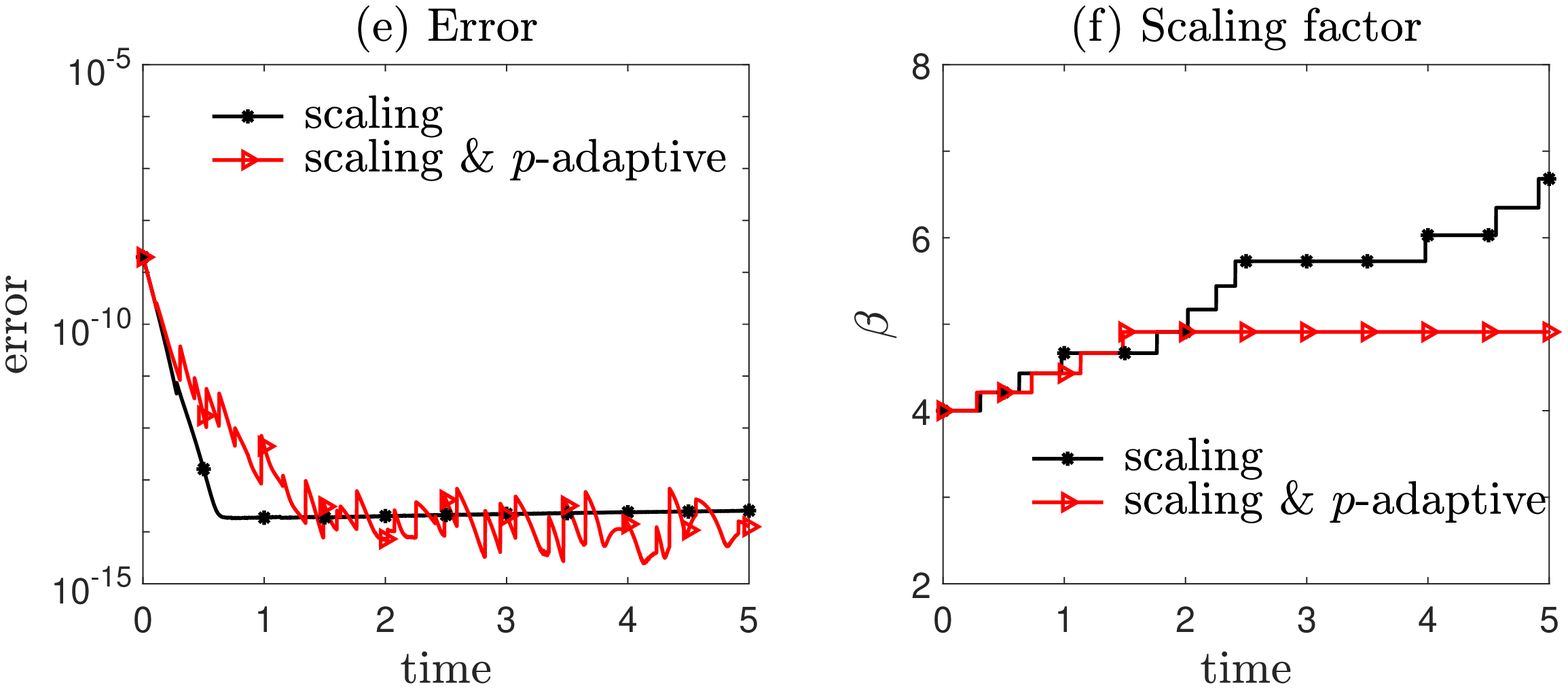}
\end{center}
\vspace{-3mm}
        \caption{\small Approximation to Eq.~\eqref{coarsescale} with
          scaling and $p$-adaptive spectral methods. Increasing
          $\beta$ by scaling can save computational burden while
          maintaining accuracy by more efficiently redistributing
          allocation points.  The approximation error is controlled
          below the initial approximation error for both scaled and
          unscaled $p$-adaptive methods, but the expansion order of
          the scaled method is smaller. On the other hand, adjusting
          the frequency indicator without decreasing $N$ will not
          achieve higher accuracy even with a much larger expansion
          order.}
     \label{fig6}
  \end{figure}

The errors and expansion orders over time are plotted in
Fig.~\ref{fig6} where the $p$-adaptive method is compared with the
non-$p$-adaptive method when scaling is applied. From
Figs.~\ref{fig6}(a) and (c) we can observe that both scaled and
unscaled methods maintain the error below the initial approximation
error. Yet, upon comparing Fig.~\ref{fig6}(b) to Fig.~\ref{fig6}(d) it
is readily seen that the scaled method leads to appropriate coarsening
while succeeding in maintaining low error, but the unscaled method
will increase $N$ when increasing the expansion order is not actually
needed, resulting in additional unnecessary computational burden. In
Fig.~\ref{fig6}(e) the scaled and $p$-adaptive spectral method with
$\eta_0=4$ is compared with the scaling-only spectral method.  We see
that the errors for both methods are almost the same but the
$p$-adaptive method can reduce unnecessary computation by decreasing
$N$ adaptively while still maintaining a low error, and the
approximation error for the $p$-adaptive method fluctuates due to a
decreasing $N$.

Fig.~\ref{fig6}(f) shows that the scaling factor $\beta$ is increased
more in the $p$-adaptive method, implying that the reason why the
$p$-adaptive method can achieve the same accuracy as non-$p$-adaptive
method with a smaller expansion order is that it can redistribute the
allocation points more efficiently.

Finally, we conclude that all
three methods: scaling, $p$-adaptive+scaling, and $p$-adaptive methods
can maintain the error well below the initial approximation error, but
the combined $p$-adaptive+scaling method can achieve this with the
smallest expansion order and is therefore the most efficient method
among them.
\end{example}

\section{Applications in solving the Schr\"{o}dinger equation}
\label{schrodinger}
In this section, we apply our adaptive spectral methods described in
Fig.~\ref{algmovingscaling} to solve the Schr\"{o}dinger equation in
unbounded domains
\begin{equation}
i\p_t \psi(x, t) = -\p_x^2\psi(x, t)+V(x)\psi(x, t) + 
V_{\rm ex}(x, t)\psi(x, t), \quad x\in\mathbb{R},
\label{Shrodinger2}
\end{equation}
which is equivalent to the PDE discussed in \cite{li2018stability}
\begin{equation}
i\p_tu(x, t) = \left[-(\p_x+iA(x, t))^2 + V(x, t)\right]u(x, t)
\label{example210}
\end{equation}
under the transformation $u(x, t) = e^{i\int_0^tV_{\rm ex}(x,
  s)\dd{s}}\psi(x, t)$.  Here, we shall use spectral methods with the
Hermite function basis. The solution is complex, so in the spectral
decomposition the coefficients of the basis functions are complex. The
major difference here is that in \cite{li2018stability} the
Schr\"{o}dinger equation is solved in a bounded domain $(x_{-}, x_+)$
with absorbing boundary conditions. Using spectral methods, we are
able to solve the Schr\"{o}dinger equation without truncating the
domain. 

We solve the weak form of Eq.~\eqref{Shrodinger2}
\begin{equation}
(\p_t \psi, v) = -i(\p_x \psi, \p_x v)+ 
((V(x) + V_{\rm ex}(x, t))\psi, v), \quad v~\in L^2(-\infty, \infty),
\end{equation}
which is to find $\Psi_{N, x_L}^{\beta}(t, x) \coloneqq
\sum_{i=0}^N\psi_{i, x_L}^{\beta}(t)\hat{\mathcal{H}}_i^{\beta}(x-x_L)$ in
$V_{N,
  x_L}^{\beta}=\textrm{span}\{\hat{\mathcal{H}}_i^{\beta}(x-x_L)\}_{i=0}^N$
satisfying the initial condition and

\begin{equation}
(\p_t \Psi_{N, x_L}^{\beta}, v) + i(\p_x \Psi_{N, x_L}^{\beta}, \p_x v)
=-i((V(x) + V_{\rm ex}(x, t))\Psi_{N, x_L}, v), \quad \forall v\in V_{N, x_L}^{\beta}.
\label{numweak2}
\end{equation} 
We denote $\boldsymbol\psi_{N, x_L}^{\beta}(t)\coloneqq
(\psi_{0, x_L}^{\beta}(t),...,\psi_{N, x_L}^{\beta}(t))$, which can be analytically
solved to advance time

\begin{equation}
\boldsymbol\psi_{N, x_L}^{\beta}(t_{n+1}) = 
\exp\left[i\int_{t_n}^{t_{n+1}}(D_N^{\beta}+V_{N, x_L}^{\beta}(t))\dd{t}\right]\boldsymbol\psi_{N, x_L}^{\beta}(t_{n}) 
\label{forwardtime}
\end{equation}
where $D^{\beta}_N\in\mathbb{R}^{(N+1)\times(N+1)}$ is a symmetric
matrix with entries
\begin{equation}
(D^{\beta}_N)_{\ell j}=\left\{
\begin{aligned}
&\beta^2\sqrt{\ell(\ell+1)} \;\;\quad\quad\quad j=\ell+2,\\
&\beta^2\sqrt{(\ell-2)(\ell-1)} \quad j=\ell-2,\\
&\beta^2\f{\ell}{2} \quad\quad\quad\quad\quad\quad\quad j=\ell,\\
&0 \quad\quad\quad\quad\quad\quad\quad\quad\; \textrm{otherwise},
\end{aligned}
\right.
\label{Ddefine}
\end{equation}
and the matrix $V_{N, x_L}^{\beta}(t)\in\mathbb{R}^{(N+1)\times(N+1)}$
has entries

\begin{equation}
(V_{N, x_L}^{\beta}(t))_{\ell j} = \int_{-\infty}^{\infty}
(V(x) + V_{\rm ex}(x, t))\hat{\mathcal{H}}_{\ell-1}^{\beta}(x-x_L)
\hat{\mathcal{H}}_{j-1}^{\beta}(x-x_L)\dd{x}.
\end{equation}
The evaluation of $\exp(i\int_{t_n}^{t_{n+1}}(D_N^{\beta}+V_{N,
  x_L}^{\beta}(t))\dd{t})\boldsymbol\psi_{N, x_L}^{\beta}(t_{n}) $ is
performed as follows. First, we denote $\tilde{V}_{N,
  x_L}^{\b}\approx\int_{t_n}^{t_{n+1}}V_{N, x_L}^{\beta}(t)\dd{t}$ where the
integration is evaluated by Gauss-Legendre formula. Therefore, when
calculating the matrix-vector product $\tilde{V}_{N, x_L}^{\b}\textbf{X}_N$ for
a vector $\textbf{X}_N\coloneqq(X_1,...,X_N)\in\mathbb{R}^{N+1}$, its
$\ell^{\rm th}$ component is

\begin{equation}
\begin{aligned}
&(\tilde{V}_{N, x_L}\textbf{X}_N)_{\ell} = \sum_{j=0}^N\sum_{s=0}^N
\hat{H}_{\ell-1}^{\b}(x_s^{\beta})\hat{H}_j^{\b}(x_s^{\beta})\bigg[
V(x_s^{\beta}+x_L) + 
\tfrac{5}{18}V_{\rm ex}(x_s^{\beta}+x_L, t_n+\tfrac{1}{2}(1-\sqrt{\tfrac{3}{5}})\dd t) \\
& \qquad \qquad \qquad + \tfrac{4}{9}V_{\rm ex}(x_s^{\beta}+x_L, t_n+\tfrac{{\rm d} t}{2})
 + \tfrac{5}{18}V_{\rm ex}(x_s^{\beta}+x_L, t_n+\tfrac{1}{2}(1+\sqrt{\tfrac{3}{5}})\dd t)\bigg]X_j \Delta{ t}
\end{aligned}
\end{equation}
where $\Delta{t}=t_{n+1}-t_n$. We can first calculate 

\begin{align}
\: & \sum_{j=0}^N\hat{H}_j^{\b}(x_s^{\beta})\bigg[V(x_s^{\beta}+x_L) 
+ \tfrac{5}{18}V_{\rm ex}(x_s^{\beta}+x_L, t_n+\tfrac{1}{2}(1-\sqrt{\tfrac{3}{5}})\dd t) \nonumber \\
\: & \qquad + \tfrac{4}{9}V_{\rm ex}(x_s^{\beta}+x_L, t_n+\tfrac{{\rm d}}{2}) 
+ \tfrac{5}{18}V_{\rm ex}(x_s^{\beta}+x_L, t_n+\tfrac{1}{2}(1+\sqrt{\tfrac{3}{5}})\dd{t})\bigg]X_j\Delta{ t}
\end{align}
for each subindex $s$; then, evaluating $(\tilde{V}_{N,
  x_L}^{\b}\textbf{X}_N)_{\ell}$ for each subindex $\ell$ will only
require an $O(N)$ operation. In this way, given any arbitrary
potentials $V(x), V_{\rm ex}(x, t)$ we can calculate $\tilde{V}_{N,
  x_L}^{\b}\textbf{X}_N$ in $O(N^2)$ operations without explicitly
calculating entries in $\tilde{V}_{N, x_L}^{\b}$. We approximate the
matrix-vector product $\exp\left[i(D_N^{\beta}\Delta{t}+\tilde{V}_{N,
    x_L})\right]\boldsymbol\psi_{N, x_L}^{\beta}(t_{n})$ in the
following way: we rewrite $\exp\left[i(D_N\Delta{t}+\tilde{V}_{N,
    x_L})\right]\boldsymbol\psi_{N, x_L}^{\beta}(t_{n}) =
\exp\left[\f{1}{m}i(D_N\Delta{t}+\tilde{V}_{N,
    x_L})\right]^m\boldsymbol\psi_{N, x_L}^{\beta}(t_{n})$, which is
introduced as the ``scaling and squaring" method in
\cite{1978Nineteen}, and approximate the matrix-vector product
$\exp\left[\f{1}{m}i(D_N\Delta{t}+\tilde{V}_{N,
    x_L})\right]\textbf{X}_N$ by truncating the infinite Taylor
expansion series
$\sum_{j=0}^{\infty}\f{1}{m^{j}j!}\left[i(D_N\Delta{t}+\tilde{V}_{N,
    x_L})\right]^j \textbf{X}_N$. Here, we take $m=6$.


As mentioned in Section~\ref{sec:intro}, two main numerical difficulties when
solving the Schr\"{o}dinger equation are the unboundedness and
oscillatory behavior of the solutions.  In fact, the solution may be
increasingly oscillatory behavior at infinity over time, making it
very hard to solve in the unbounded domain. However, with our adaptive
spectral methods, we can efficiently solve the Schr\"{o}dinger
equation in unbounded domains accurately and capture these
oscillations.

We shall revisit the two numerical examples discussed in
\cite{li2018stability}. In the following examples, curves labeled
``adaptive'' indicate that scaling, moving, and $p$-adaptive
techniques are all applied as in the algorithm described in
Fig.~\ref{algmovingscaling}, while curves labeled ``$p$-adaptive''
mean that Alg.~\ref{algrefining} is used without scaling or moving;
similarly, curves labeled ``scaling \& $p$-adaptive'' are obtained
when the $p$-adaptive and scaling techniques are applied.  Curves
labeled ``scaling \& moving'' or ``scaling" correspond to applying the
scaling and moving techniques or the scaling technique without the
$p$-adaptive subroutine.  The ``non-adaptive'' curves are obtained
when we do not apply any of the scaling, moving, or $p$-adaptive
techniques.

\begin{figure}
\begin{center}
      \includegraphics[width=5in]{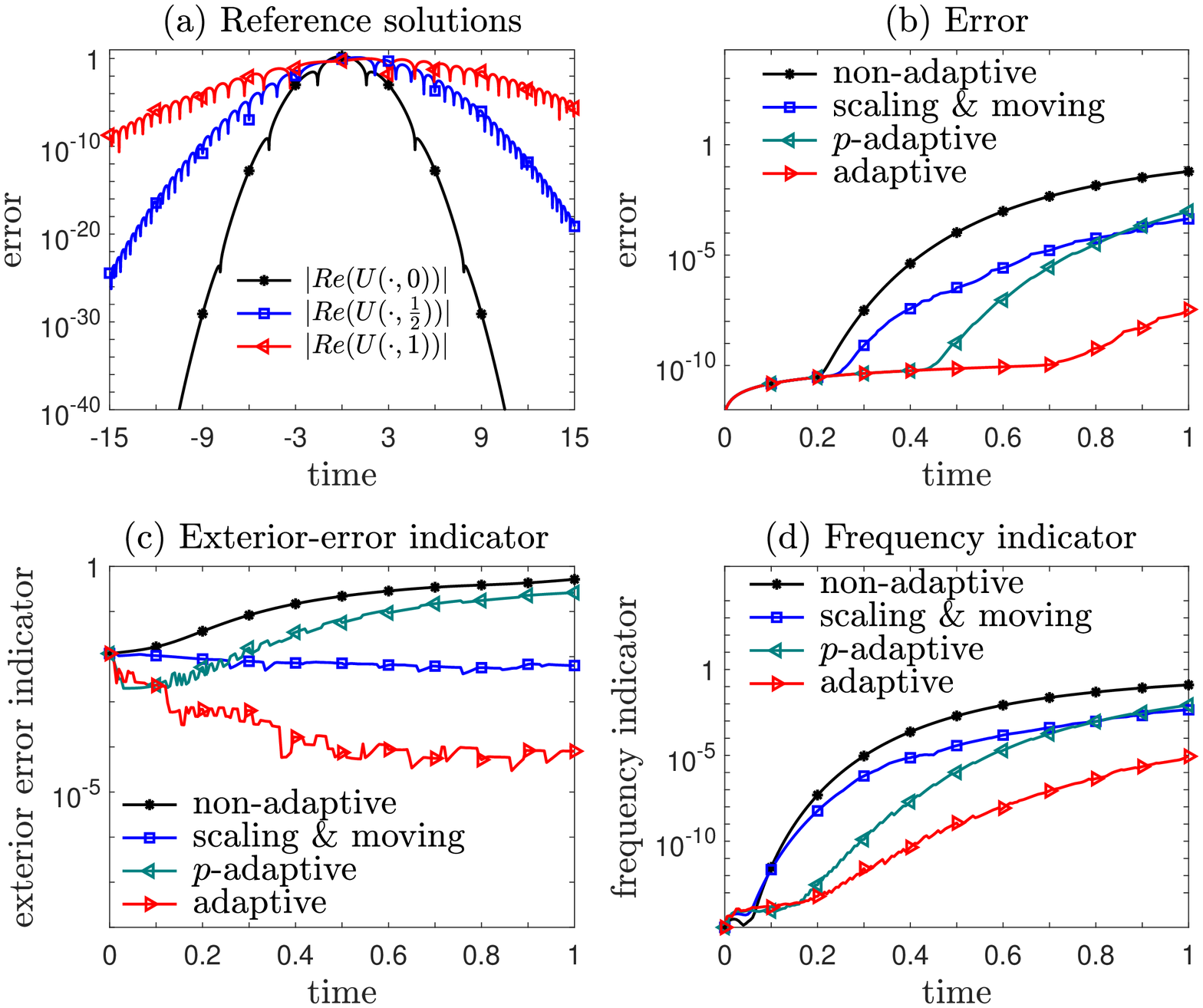}
\includegraphics[width=5in]{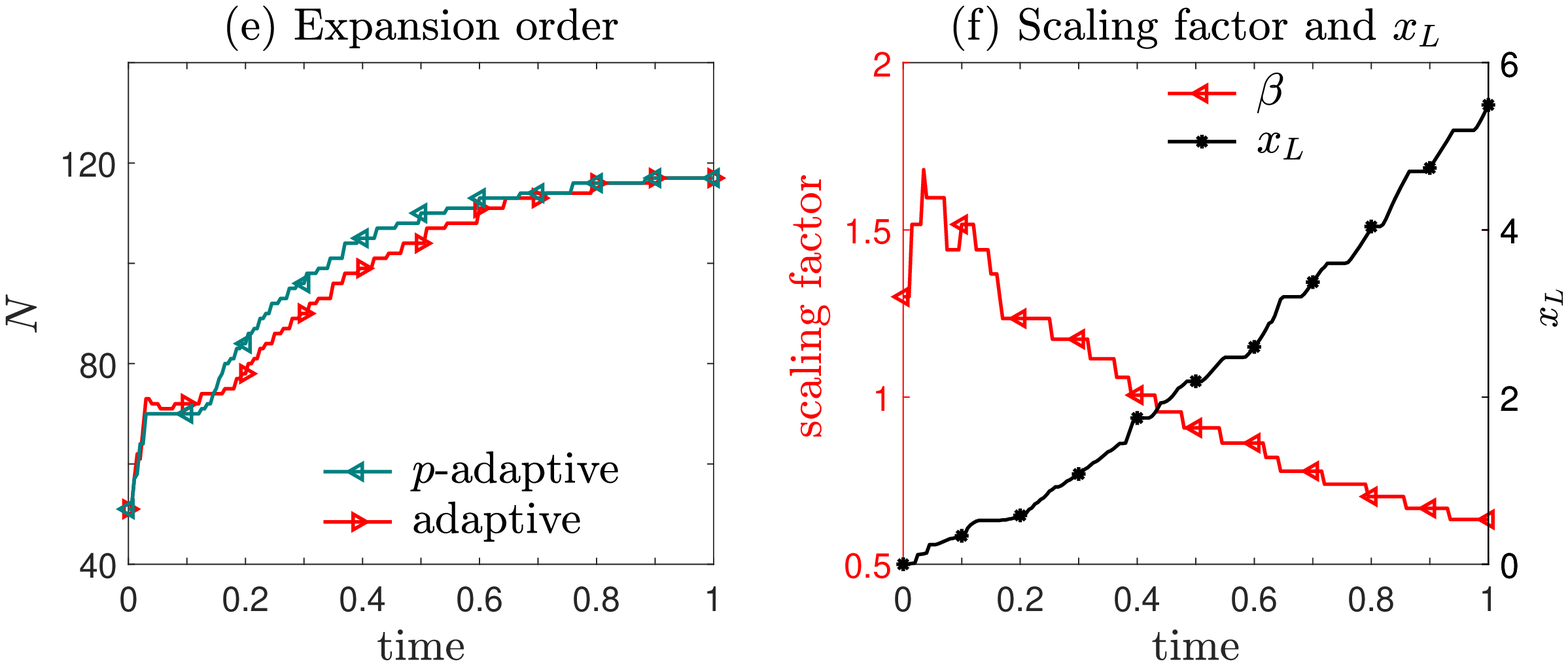}
\end{center}
\vspace{-3mm}
        \caption{\small Numerically solving the Schr\"{o}dinger
          equation with vanishing potentials.  Applying scaling, moving, and
          $p$-adaptive techniques can successfully capture diffusive,
          advective, and oscillatory behavior of the solution and
          yields an accurate numerical solution that prevents the
          frequency indicator from growing too fast.  The
          exterior-error indicator is also kept small by moving the
          basis functions rightward to avoid a deteriorating
          approximation at $\infty$. Failure to incorporate any of the
          moving, scaling, or $p$-adaptive techniques results in a
          much larger error.}
     \label{fig4}
\end{figure}

\begin{example}
\label{ex:simulatebounded}
\rm We numerically solve the Schr\"{o}dinger equation which is solved
in Example 1 of \cite{li2018stability} and take $V=V_{\rm ex}=0$ in
Eq.~\eqref{Shrodinger2}, admitting the analytic solution
\begin{equation} 
\Psi(x, t) = \frac{1}{\sqrt{\zeta+it}}\exp\left[ik(x-kt)-
\frac{(x-2kt)^2}{4(\zeta+it)}\right],
\end{equation}
where $k$ is related to the propagation speed of the beam and $\zeta$
determines the width of the beam.  The absolute values of the real
part of $\Psi(x, t=0,0.5,1)$ are plotted in Fig.~\ref{fig4}(a),
illustrating the increasingly oscillatory and diffusive behavior in
the rightward propagating solution. Treatment of this solution will
thus require scaling, moving, and $p$-adaptive techniques.  The
imaginary parts of the reference solution (not plotted) over time are
also increasingly oscillatory.  We shall apply the algorithm described
in Fig.~\ref{algmovingscaling}. We set $\zeta=0.3, k=1$, and
initialize $N=50$ at $t=0$. Other parameters are set to $q =
\nu^{-1}=0.95, \mu=1.0002, d_0=0.005, \underline{\beta}=0.3,
\overline{\beta}=2, d_{\max}=0.1, N_{\max}=6, N_{\min}=0, \eta=1.1,
\gamma=1.05$, and $\Delta{t}=0.005$. Note that with zero potential,
Eq.~\eqref{forwardtime} reduces to

\begin{equation}
\boldsymbol\psi_N^{\beta}(t_{n+1}) =
\exp(iD_N^{\beta}\dd{t})\boldsymbol\psi_N^{\beta}(t_{n}).
\label{potential0}
\end{equation}
When all four techniques are applied, the error is the smallest (shown
in Fig.~\ref{fig4}(b)) since we can keep the exterior error indicator
in $(x_R, \infty)$ small (shown in Fig.~\ref{fig4}(c)) by matching the
solution's intrinsic advection. We can simultaneously prevent the
frequency indicator from growing too fast (shown in
Fig.~\ref{fig4}(d)), thus ensuring a small error bound.

From the reference solution it can be observed that increasing the
expansion order over time is an intrinsic requirement and failure to
do so prevents the capture of the increasing oscillations, leading to
a huge error.  As the function becomes increasingly oscillatory as
$x\rightarrow\infty$, moving the basis rightward requires
correspondingly more refinement (shown in Fig.~\ref{fig4}(e)).
However, the $p$-adaptive method alone cannot compensate for the
inability to capture diffusion and advection, resulting in an
inaccurate approximation. We have also checked that apart from
what is shown in Fig.~\ref{fig4}, applying any single scaling, moving
or $p$-adaptive technique, or combining any two of them will all
result in a much larger error than employing all three techniques
indicated in Fig.~\ref{algmovingscaling}.

\end{example}

\noindent Finally, we numerically solve the Schr\"{o}dinger equation
with non-vanishing potentials. %
\begin{example}
\label{ex:schrodinger}
\rm We numerically solve the following standard Schr\"{o}dinger equation Eq.~\eqref{Shrodinger2} equivalent to Example 2 in
\cite{li2018stability} with potentials
\begin{equation}
V_{\rm ex}(x, t) = \f{50}{\sqrt{\pi}}\sin(10t)\int_{-\infty}^x \!\!\exp(-z^2)\dd{z}, 
\,\, V(x) = -10\left[e^{-10(x-1)^2}+e^{-10(x+1)^2}\right].
\end{equation}
Given an even function as the initial condition for Example 2 in
\cite{li2018stability}, the solution is also an even function and the
solution of Eq.~\eqref{Shrodinger2} obeys $|\psi(-x, t)|=|\psi(x,
t)|$. No bias towards $-\infty$ or $+\infty$ is preferred.  Therefore,
we use the Hermite function basis and apply the algorithm described in
Fig.~\ref{algmovingscaling} but deactivate the moving technique by
setting $d_{\rm max}=0$.

We set the initial condition to be the same as that of Example
\ref{ex:simulatebounded} and set $\eta=1.025, \gamma=1, q=0.95,
\nu=q^{-1}, N_{\min}=0, N=200, \underline{\beta}=0.3,
\overline{\beta}=2$ and $\beta_0=1.3$ at $t=0$ with the maximal
expansion order increment for each step $N_{\max}=20$.
  \begin{figure}[h!]
\begin{center}
      \includegraphics[width=5.2in]{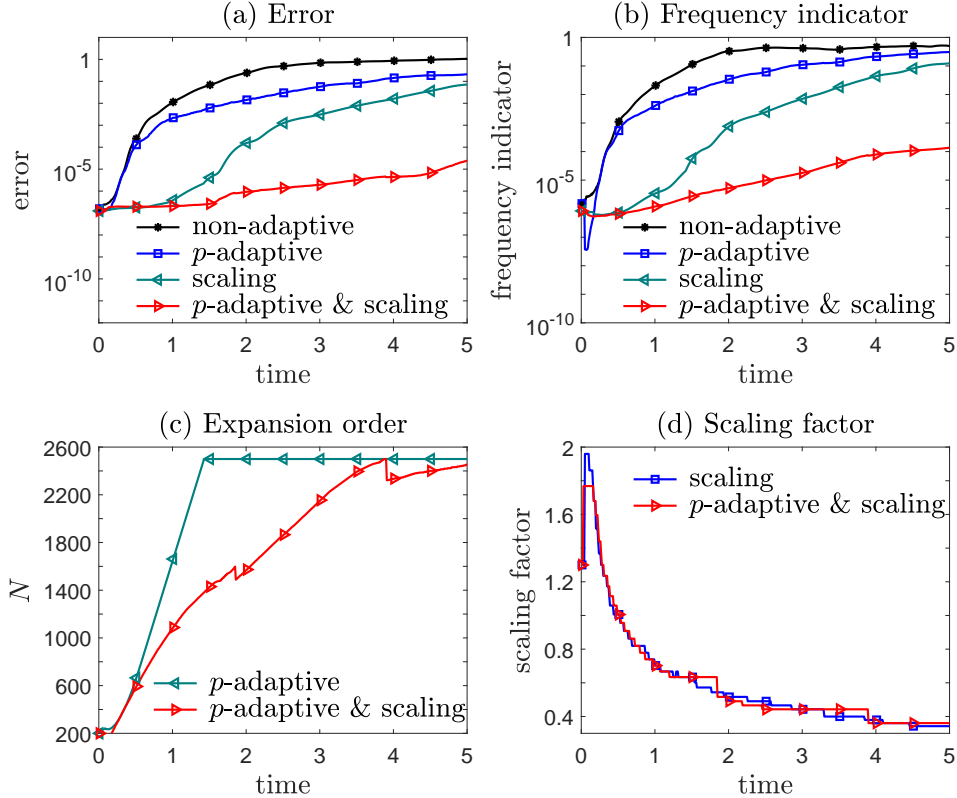}
\end{center}
\vspace{-3mm}
        \caption{\small Numerically solving the Schr\"{o}dinger
          equation with non-vanishing potentials. Rapidly increasing
          oscillations of the solution over time requires much
          refinement and proper scaling to maintain accuracy. It is
          again verified that proper scaling can avoid unnecessary
          refinement and avoid unnecessary computational burden by
          adaptively adjusting the scaling factor. Without scaling,
          the expansion order soon reaches the upper bound for $N$
          (the expansion order of the reference solution) and the
          approximation soon deteriorates due to an inability to
          further increase $N$ or adjust $\beta$ and maintain a low
          frequency indicator.  Failure to accommodate the
          $p$-adaptive technique will also result in a larger error
          because of an inability to capture the oscillatory
          behavior.}
     \label{fig5}
\end{figure}

The reason why we set $\gamma=1$ is that the expansion order $N$ needs
to be increased quickly to catch up with the highly increasingly
oscillatory behavior of the numerical solution. We set a uniform
timestep $\Delta{t}=0.01$. We use the numerical solution solved with a
fixed $N=2500$ and only the scaling technique activated as the
reference solution.  For the $p$-adaptive method, we have added an
additional restriction that the expansion order cannot surpass the
expansion order of the reference solution $N=2500$.

We can easily see that the spectral method with both scaling and
$p$-adaptive techniques outperforms the non-adaptive spectral method
or with only one of these two techniques employed (shown in
Fig.~\ref{fig5}(a)).  The frequency indicator of using both 
scaling and $p$-adaptive techniques is also the smallest
(Fig.~\ref{fig5}(b)), and the similarity between the frequency
indicator and error is again confirmed as stated in \cite{Xia2020b}.
Moreover, the unscaled method will result in a larger expansion order
(shown in Fig.~\ref{fig5}(a)), leading to excessive refinement with no
improvement in accuracy (shown in Fig.~\ref{fig5}(a)). In this
example, the coarsening procedure will not lead to a large increase of
frequency indicator and does not significantly compromise accuracy
(shown in Figs.~\ref{fig5}(b) and (c)).  Finally, the scaling factors
of the $p$-adaptive spectral method and the non-$p$-adaptive spectral
method trend similarly over time; they both decrease after
experiencing an initial, transient increase (Fig.~\ref{fig5}(d)).
\end{example}

\section{Summary and Conclusion}
\label{conclusion}
In this paper, we proposed a frequency-dependent $p$-adaptive
technique that adjusts the expansion order for spectral methods. We
demonstrated its applicability to time-dependent problems with varying
oscillatory behavior.  In order to develop efficient numerical methods
for problems requiring solutions in unbounded domains, we also
combined the $p$-adaptive technique with scaling ($r$-adaptivity) and
moving ($h$-adaptivity) methods to devise a complete adaptive spectral
method that can successfully deal with diffusion, advection, and
oscillation.

The relationship between scaling and $p$-adaptive techniques for
spectral methods in unbounded domains, both of which depend on the
same frequency indicator, is also investigated.  We successfully applied
our adaptive spectral method to numerically solve Schr\"odinger's
equation.  The associated solutions are highly oscillatory in the
whole domain, posing numerical difficulties for existing numerical
methods that truncate the domain.

For future research, the relationship among the adaptive techniques
for spectral methods, scaling, moving, refinement and coarsening, can
be further studied and rigorous numerical analysis for these
techniques should be investigated.  Furthermore, fast algorithms with
mapped Chebyshev polynomials for solving PDEs in unbounded domains
have been developed using the fast Fourier transform
\cite{sheng2020fast}, but there lacks fast and efficient algorithms
exploiting Laguerre and Hermite basis functions, particularly for
higher-dimensional problems. Thus, generalizing these adaptive methods
for mapped Jacobi polynomials may be a compelling future research
direction.

\section*{Acknowledgments}

MX and TC acknowledge support from the National Science Foundation
through grant DMS-1814364 and the Army Research Office through grant
W911NF-18-1-0345.  SS acknowledges financial support from the National
Natural Science Foundation of China (Nos.~11822102, 11421101) and
Beijing Academy of Artificial Intelligence (BAAI).  Computational
resources were provided by the High-performance Computing Platform at
Peking University.

\bibliographystyle{elsarticle-num} 

\end{document}